\documentclass{conm-p-l}

\newtheorem{theorem}{Theorem}[section]
\newtheorem{proposition}[theorem]{Proposition}

\newtheorem{lemma}[theorem]{Lemma}

\theoremstyle{definition}
\newtheorem{definition}[theorem]{Definition}
\newtheorem{example}[theorem]{Example}

\newtheorem{problem}[theorem]{Problem}

\theoremstyle{remark}

\numberwithin{equation}{section}

\newcommand{\x}{\times}
\newcommand{\cH}{\mathcal H}
\newcommand{\cB}{\mathcal B}
\newcommand{\cN}{\mathcal N}
\newcommand{\cM}{\mathcal M}
\newcommand{\bbR}{\mathbb R}
\newcommand{\bbC}{\mathbb C}
\newcommand{\bbT}{\mathbb T}
\newcommand{\bbZ}{\mathbb Z}
\newcommand{\fg}{\mathfrak g}

\begin{document}

\title{Compact Quantum Metric Spaces} 
\author{Marc A. Rieffel}
\address{Department of Mathematics, University of California
Berkeley, CA 94720-3840}
\email{rieffel@math.berkeley.edu}
\thanks{The author's research was supported in part by NSF Grant DMS-0200591.}
\date{June 30, 2003}
 \copyrightinfo{2003}
    {American Mathematical Society}


\begin{abstract}
We give a brief survey of many of the high-lights of our present understanding
of the young subject of quantum metric spaces, and of quantum
Gromov-Hausdorff distance between them. We include examples.
\end{abstract}

\maketitle

My interest in developing the theory of compact quantum metric spaces
was stimulated by certain statements in the high-energy physics and string-theory
literature, concerning non-commutative spaces that converged to other
(possibly noncommutative) spaces. These statements appeared to me to deserve a
more precise formulation.
(See the references in the introductions of \cite{R6} \cite{R7}.) Here
I will just give a brief survey of some of the main developments in this very young
subject. I will also indicate some of the main classes of examples which have 
been explored so far. I include a few related arguments which are not quite 
in place in the existing literature. 

\section{The definition of compact quantum metric spaces}

The concept of a quantum metric space has its origins in Connes' 
paper \cite{C1}
of 1989 in which he first proposes using Dirac operators as the vehicle for
metric data in non-commutative geometry.  He was motivated by his observation
that for a compact spin Riemannian  manifold one can recover its smooth
structure, its Riemannian metric, and much else, directly from its 
standard Dirac
operator.  This led him to the concept of a spectral triple, $(A,\cH,D)$,
consisting of a $*$-algebra $A$ represented by bounded operators on a Hilbert
space $\cH$, and of a (usually unbounded) self-adjoint operator, $D$, on $\cH$
such that the commutator $[D,a]$ is a bounded operator for each $a \in A$.
Connes also requires that $D$ have compact resolvant.  (Spectral 
triples are very
closely related to ``unbounded $K$-cycles'' or ``unbounded Fredholm modules'',
the difference being that for spectral triples the representation of 
$A$ must be
faithful.)  Connes pointed out that from a spectral triple one obtains a metric,
$\rho_D$, on the state space, $S(A)$, of $A$ by means of the formula
\[
\rho_D(\mu,\nu) = \sup\{|\mu(a) - \nu(a)|: \|[D,a]\| \le 1\},
\]
where the value $+\infty$ may occur.  But he did not explore this 
metric much for
non-commutative algebras.

In \cite{R4} \cite{R5} I pointed out that, motivated by what happens for
ordinary compact metric spaces, it is natural to desire that the topology on
$S(A)$ determined by $\rho_D$ coincide with the weak-$*$ topology.  To be more
specific about this motivation, let $(X,\rho)$ be an ordinary compact metric
space.  It is common to define a seminorm, $L_{\rho}$, the Lipschitz 
seminorm, on
the algebra $A = C(X)$ of continuous functions on $A$, by
\[
L_{\rho}(f) = \sup\{|f(x)-f(y)|/\rho(x,y): x \ne y\},
\]
where the value $+\infty$ is permitted.  The metric $\rho$ can be 
recovered from
$L_{\rho}$ by
\[
\rho(x,y) = \sup\{|f(x)-f(y)|: L_{\rho}(f) \le 1\}.
\]
Even more, Kantorovich \cite{Kn} \cite{KR} used $L_{\rho}$ to define 
a metric on
$S(A)$, which now is identified with the space of all probability measures on
$X$, by
\[
\rho_L(\mu,\nu) = \sup\{|\mu(f)-\nu(f)|: L_{\rho}(f) \le 1\}.
\]
It is clear that $\rho_L$ extends $\rho$ from the set of point-measures to the
set of all probability measures.  Kantorovich showed, among other properties,
that the topology on $S(A)$ from $\rho_L$ coincides with the weak-$*$ topology.

Since the data $(C(X),L_{\rho})$ is equivalent to the data $(X,\rho)$, we see
that we have reformulated the notion of a metric in terms of just the 
commutative
$C^*$-algebra $C(X)$, by means of $L_{\rho}$.  It is then 
natural to try to
formulate metric data for a non-commutative unital $C^*$-algebra, $A$, by means
of a suitable seminorm on $A$ which should play the role of $L_{\rho}$. For
Connes' spectral triples the role of $L_\rho$ is played, in effect, by
$L(a) = \|[D,a]\|$.

However, for technical flexibility it is useful to consider the more general
situation in which $A$ is just an order-unit space.  The definition of an
order-unit space is due to Kadison \cite{Kd}, influenced by work of Stone.  For
our present purposes it is sufficient to know that any order-unit space can be
realized concretely (in many ways) as a real linear space of self-adjoint
operators on a Hilbert space, containing the identity operator (the order unit), and 
equipped with
the usual partial ordering and norm on operators.  An order-unit 
space $A$ has a
state-space, $S(A)$, defined just as for $C^*$-algebras.  We can consider a
seminorm, $L$, on an order-unit space, and use it to define a metric, $\rho_L$,
on $S(A)$, much as above, by
\[
\rho_L(\mu,\nu) = \sup\{|\mu(a)-\nu(a)|: L(a) \le 1\}.
\]

To be useful we need some conditions on $L$.  As first condition, we 
require for
convenience that $L(1) = 0$ where $1$ is the order-unit.  This is 
well-motivated
by the case of ordinary metric spaces.  But the main requirement 
which we make is
that the topology on $S(A)$ from $\rho_L$ should coincide with the weak-$*$
topology.  This is motivated by the result of Kantorovich for ordinary metric
spaces which we mentioned above.  We view the role of this requirement as
follows.  Given a compact topological space $X$, there are plenty of metrics on
$X$ as a set, but when we say ``a metric on $X$'' we nearly always also have in
mind ``whose topology coincides with the given topology on $X$''.  We consider
the requirement that $\rho_L$ give $S(A)$ the weak-$*$ topology to be 
the analog
of this idea for seminorms on order-unit spaces.

We can always restrict to the subspace on which $L$ takes finite values.  This
does not change $\rho_L$.  (We do not require our order-unit spaces to be
complete.)  Thus we make \cite{R6}:

\begin{definition}
Let $A$ be an order-unit space, and let $L$ be a seminorm on $A$ taking finite
values.  We say that $L$ is a {\em Lip-norm} if
\begin{itemize}
\item[1)] $L(1) = 0$.
\item[2)] The topology on $S(A)$ from $\rho_L$ coincides with the weak-$*$
topology.
\end{itemize}
By a {\em compact quantum metric space} we mean a pair $(A,L)$ where $A$ is an
order-unit space and $L$ is a Lip-norm on $A$.
\end{definition}

It is easily seen that condition $2$ implies that if $L(a) = 0$ then 
$a \in \bbR
1$.  If $L$ is actually a seminorm on a dense subalgebra of a unital
$C^*$-algebra, $A$,
or more generally on a (complex) self-adjoint linear subspace of bounded
operators, then we require that $L(a) = L(a^*)$ for all $a \in A$. 
Then a simple
argument (in section 2 of \cite{R6})
shows that $L$, and the restriction of $L$ to the order-unit space $A_{sa}$ of
self-adjoint elements of $A$, determine the same metric on $S(A)$.  We can then
require that either one gives $S(A)$ the weak-$*$ topology.

Aside from the motivation given by the case of ordinary metric spaces, it is
reasonable to ask why we care that $\rho_L$ give $S(A)$ the weak-$*$ topology.
The topic of quantum metric spaces is still in its infancy, and it is far too
soon to know what demands new important examples will bring.  At the moment the
main answer to this question which I have is that this requirement permits an
effective notion of Gromov--Hausdorff distance for compact quantum metric
spaces.  We will discuss this in Section~5, followed by examples in Section~6.

\section{Sources of examples}

Of course, it is condition $2$ in the definition of a Lip-norm which may be
difficult to verify in naturally occurring examples.  It can be reformulated in
terms of $A$ itself in the following way \cite{R4}.  Let ${\tilde A} 
= A/\bbC 1$,
equipped with the corresponding quotient norm $\|\cdot\|^{\sim}$.

\begin{theorem}
Let $L$ be a seminorm on an order-unit space $A$ such that $L(1) = 0$.  Set
$\cB_1 = \{a \in A: L(a) \le 1\}$.
\begin{itemize}
\item[a)] Then $\rho_L$ gives $S(A)$ finite diameter iff the image of
$\cB_1$ in ${\tilde A}$ is bounded.
\item[b)] And $\rho_L$ gives $S(A)$ the weak-$*$ topology iff the image of
$\cB_1$ in ${\tilde A}$ is totally bounded (for $\|\cdot\|^{\sim}$).
\end{itemize}
\end{theorem}

It can be shown by a somewhat unnatural construction \cite{R18} that every 
separable
order-unit space has an abundance of Lip-norms (finite on dense subspaces).
Certain other constructions are also known \cite{R18} \cite{Ckb2} 
\cite{AnC}.  But most
interesting constructions so far have come from actions of groups.  In fact, we
will eventually see that in a certain sense all Lip-norms come from actions of
Lie groups, and in fact actions of $\bbR$.

Let $A$ be a unital $C^*$-algebra, let $G$ be a locally compact group, and let
$\alpha$ be a strongly continuous action of $G$ on $A$ by automorphisms.  To
define a Lip-norm we must somewhere put in metric data.  We do this 
by choosing a
continuous length-function $\ell$ on $G$.  This means that $\ell$ has values in
$\bbR^+$, that $\ell(xy) \leq \ell(x) + \ell(y)$, that $\ell(x^{-1}) = 
\ell(x)$, and
that $\ell(x) = 0$ iff $x = e$.  Then we can define a seminorm, $L$, on $A$ by
\[
L(a) = \sup\{\|\alpha_x(a)-a\|/\ell(x): x \ne e\}.
\]
It can be seen by a ``smoothing argument'' that $\{a: L(a) < 
\infty\}$ is a dense
$*$-subalgebra of $A$ if $\ell$ is a proper function.  
A proof of this for compact $G$ is
given in proposition $2.2$ of \cite{R4}.  But that proof can be easily modified
to apply also to locally compact groups if $\ell$ is proper, as 
follows.  For each
$n$ let $g_n = (n^{-1} - \ell)^+ \in C_c(G)$.  Let $f_n = c_ng_n$ 
where $c_n$ is
the positive constant such that $\|f_n\|_1 = 1$.  Then $\{f_n\}$ is an approximate
identity for $L^1(G)$, consisting of functions which are Lipschitz for the
right-invariant metric on $G$ defined by $\rho(x,y) = \ell(xy^{-1})$.  Simple
calculations then show that $L(\alpha_{f_n}(a)) < \infty$ for any $a \in A$,
where $\alpha_{f_n}(a) = \int f_n(x)\alpha_x(a)dx$.

If we are to have that $L(a) = 0$ only when $a \in \bbC 1$, it is clear that
$\alpha$ must be {\em ergodic}, in the sense that the only elements 
of $A$ which
are invariant under $\alpha$ are those in $\bbC 1$.  It is shown in \cite{R4}
that:

\begin{theorem}
If $G$ is compact and $\alpha$ is ergodic, then $L$ (restricted to $A_{sa}$)
is a Lip-norm.
\end{theorem}

This easily fails if $G$ is not compact.

If $G$ is a Lie group, then we can consider the set $A^{\infty}$ of smooth
elements of $A$ for the action $\alpha$, that is, elements $a \in A$ such that
the function $x \mapsto \alpha_x(a)$ on $G$ is infinitely differentiable.  Then
$A^{\infty}$ is a dense $*$-subalgebra of $A$ \cite{BrR}, and $\alpha$ gives a
representation (also denoted by $\alpha$) of the Lie algebra, $\fg$, of $G$ by
derivations on $A^{\infty}$.  Given $a \in A^{\infty}$, we can define its
differential, $da$, to be the operator from $\fg$ into $A^{\infty}$ defined by
$da(X) = \alpha_X(a)$ for $X \in \fg$.  If we bring metric data into the
situation by choosing a norm on $\fg$, then the norm of $da$ as an 
operator from
$\fg$ to $A$ is defined, and we can define a seminorm $L$ on 
$A^{\infty}$ by $L(a)
= \|da\|$.  If $G$ is compact then we can again show \cite{R4} that $L$ is a
Lip-norm, by reducing this situation to that of Theorem $2.2$.

In general, compact groups have many ergodic actions on unital $C^*$-algebras.
Here are three constructions.
\begin{itemize}
\item[1)] If $U$ is an irreducible unitary representation of a 
compact group $G$
on a Hilbert space $\cH$, then we can define an action, $\alpha$, of $G$ on $A =
\cB(\cH)$ by $\alpha_x(T) = U_xTU_x^*$.  This action is ergodic.
\item[2)] If $G$ is a compact Abelian group and if $c$ is a $2$-cocycle with
values in $\bbT$ on its (discrete) dual group ${\hat G}$, then we can define
\cite{ZM} the twisted group $C^*$-algebra $A = C^*({\hat G},c)$.  There is a
natural action, $\alpha$, of $G$ on $A$ given by $\alpha_x(f)(\gamma) = \langle
x,\gamma\rangle f(\gamma)$ for $f \in \boldsymbol{\ell}^1({\hat G})$ 
and $\gamma \in {\hat
G}$.  This action (called the ``dual action'') is ergodic.
\item[3)] If $H$ is a closed subgroup of $G$ and if $\beta$ is an action of $H$
on a unital $C^*$-algebra $B$, then we can form the {\em induced} $C^*$-algebra
$A$, consisting of the continuous $B$-valued functions $F$ on $G$ which satisfy
the condition that $F(xs) = \beta_s(F(x))$ for $x \in G$ and $s \in H$.  The
action $\alpha$ of left translation by elements of $G$ carries $A$ into itself,
and it is easily verified that if $\beta$ is ergodic then so is $\alpha$.
\end{itemize}
We can then combine this inducing construction with the two previous
constructions to obtain many ergodic actions.  But I do not know of other
constructions of ergodic actions of compact groups.  
For $G$ Abelian the second and third
constructions give all ergodic actions \cite{OPT}.  But for non-Abelian compact
groups there is no known classification of the possible ergodic actions.  For
example, it seems still to be unknown whether $SU(4)$ has
any ergodic actions other than those given by the above constructions \cite{Wa}
\cite{Wa2} \cite{Wa3}, and in particular, any ergodic
actions on infinite-dimensional unital
$C^*$-algebras which are simple.

Anyway, Theorem $2.2$ combined with the above three constructions gives many
examples of compact quantum metric spaces.  We note in particular that for the
non-commutative tori, $A_{\theta}$, which come from the second construction
described above, the dual action is an action of $\bbT^d$.  Thus any continuous
length function on $\bbT^d$ (of which there is an abundance) gives Lip-norms on the $A_{\theta}$'s. 
Here $\theta$
is a real $d \x d$ skew-symmetric matrix, $c_{\theta}$ is the bicharacter on
$\bbZ^d$ defined by $c_{\theta}(x,y) = e^{2\pi i(x \cdot \theta y)}$, and
$A_{\theta}$ is \cite{R3} the twisted group $C^*$-algebra
$C^*(\bbZ^d,c_{\theta})$.

When $G$ is not compact, one will need stronger conditions on $\alpha$ in order
to obtain a Lip-norm, but this direction has barely been explored.  I know of
only two relevant papers, both of which deal with the action of the Heisenberg
Lie group on the non-commutative Heisenberg manifolds \cite{R8}.  In the first,
by Weaver \cite{W3}, the metric structures come from sub-Riemannian 
metrics, and
so it does not quite fit into the above framework.  The second, by Chakraborty
\cite{Ckb}, uses a Lip-norm which is not defined in terms of the operator norm,
so again does not fit exactly into the above framework.  Thus it is 
still unclear
what happens for non-commutative Heisenberg manifolds within our present
framework.

We find further examples of compact quantum metric spaces by a
simple application of Theorem
$2.1$b to obtain:

\begin{proposition}  Let $L$ be a Lip-norm on an order-unit space
$A$, and let $B$ be a subspace of $A$ which contains the order unit 
(so that $B$
is an order-unit space). Then the restriction of $L$ to $B$ is a Lip-norm.
\end{proposition}

In particular, if $A$ is (the self-adjoint part of) a dense $*$-subalgebra of a
unital $C^*$-algebra, and if $B$ is a unital $*$-subalgebra of $A$, then any
Lip-norm on $A$ restricts to a Lip-norm on $B$.  When used in conjunction with
Theorem $2.2$ this gives many more examples of compact quantum metric 
spaces.  We
will only describe one specific class of examples here.  Let $A_{\theta}$ be the
non-commutative $2$-torus with unitary generators $U$ and $V$ satisfying the
relation $VU = e^{2\pi i\theta}UV$ for some real number $\theta$.  As 
mentioned above,
there is a natural ergodic action $\alpha$ of $\bbT^2$ on $A_{\theta}$, the 
dual action.
For any continuous length function on $\bbT^2$ we obtain a Lip-norm $L$ on
$A_{\theta}$ by using $\alpha$.  The algebraic $*$-subalgebra, $A_{\theta}^f$,
generated by $U$ and $V$ is carried into itself by $\alpha$, and $L$ 
is finite on
$A_{\theta}^f$.  There is a (unique) involutory automorphism, $\beta$, of
$A_{\theta}$ determined by $\beta(U) = U^{-1}$ and $\beta(V) = 
V^{-1}$.  Thus the
$2$-element group acts on $A_{\theta}$, and carries $A_{\theta}^f$ into itself.
Let $B_{\theta}$ be the fixed point algebra of this automorphism.  Then
$B_{\theta}$ is not carried into itself by $\alpha$, so that we can not apply
Theorem $2.2$ directly to $B_{\theta}$.  But $B_{\theta} \cap A_{\theta}^f$ is
easily seen to be a dense $*$-subalgebra of $B_{\theta}$, and the 
restriction of
$L$ to $B_{\theta} \cap A_{\theta}^f$ will be a Lip-norm by
Proposition 2.3.  Now the 
$B_{\theta}$'s
were extensively studied by Bratteli, Elliott, Evans and Kishimoto, under the
name ``non-commutative spheres''.  Bratteli and Kishimoto showed 
\cite{BrK} that
$B_{\theta}$ is actually an AF $C^*$-algebra when $\theta$ is 
irrational.  So we
have obtained in this way some fairly natural and interesting examples of
Lip-norms on AF $C^*$-algebras. 

Another interesting class of examples of
Lip-norms on AF $C^*$-algebras
is studied in \cite{Kr2}, in connection with the development of a notion
of ``dimension'' for compact quantum metric spaces, and a notion of 
entropy for automorphisms of ($C^*$-algebraic) compact quantum metric spaces.
Also among the examples to which this notion of entropy is applied
are automorphisms of non-commutative tori.

\section{Dirac operators}

Let $G$ be a Lie group, not necessarily compact, and let $\beta$ be 
an action of
$G$ on a unital $C^*$-algebra $B$.  For a given continuous length 
function on $G$
we can define a seminorm, $L$, on $B$ as before.  But now we will not 
assume that
$\beta$ is ergodic, and so we may have a large subalgebra of $B$ on 
which $L$ takes
value $0$.  It may nevertheless happen that for suitable unital $*$-subalgebras
$A$ (or even order-unit subspaces) of $B$, which are {\em not} carried into
themselves by $\beta$, the restriction of $L$ to $A$ is a Lip-norm.

Actually, up to now the only situation that I know of in which this possibility
has been used is that in which we have a (strongly-continuous) unitary
representation $U$ of $G$ on a possibly infinite-dimensional Hilbert space
$\cH$.  Then, much as in the first construction of Section~2 we 
define an action,
$\beta$, of $G$ on $\cB(\cH)$ by $\beta_x(T) = U_xTU_x^*$.  In general this
action will not be strongly continuous.  We let $B\ (= B_{\beta})$ be the largest
subalgebra of $\cB(\cH)$ on which $\beta$ is strongly continuous.  (Then $B$ is
weak-operator dense in $\cB(\cH)$ by lemma $7.5.1$ of \cite{Pdr}.) 
Furthermore,
$B$ is carried into itself by $\beta$.  Thus we can use $(B,\beta)$ as in the
previous paragraph, and ask, for any given unital $C^*$-subalgebra 
$A$ of $B$ (or
any order-unit subspace of $B$) whether the restriction to $A$ of a seminorm on
$B$ coming from a continuous length function on $G$ is a Lip-norm.

But actually, up to now, the only situation that I know of in which this
possibility has been used is that in which $G = \bbR$, with its usual length
function.  This is the very important case of Connes' Dirac operators.  Any
strongly continuous unitary representation $U$ of $\bbR$ has an infinitesimal
generator, which we will denote here by $D$, and which is a (usually unbounded)
self-adjoint operator on $\cH$.  By means of the functional calculus for
unbounded self-adjoint operators we have $U_t = e^{itD}$ for $t \in 
\bbR$.  Let
$\beta$ be defined in terms of $U$ as above, and let $B_{\beta}$ be the
$C^*$-subalgebra of $\cB(\cH)$ on which $\beta$ is strongly continuous.  Let $T
\in B_{\beta}$, and assume further that $t \mapsto \beta_t(T)$ is once
differentiable for the operator norm.  Just from the definition of 
$D$ being the
infinitesimal generator of $U$ it follows easily that $T$ carries the domain of
$D$ into itself and that $[D,T]$ is bounded on that domain and so extends to a
bounded operator on $\cH$.  Conversely, if $T$ carries the domain of $D$ into
itself and if $[D,T]$ is bounded, then $T$ is in $B_{\beta}$. (See the 
first line
of the proof of corollary $10.16$ of \cite{GVF}.)

If $t \mapsto \beta_t(T)$ is twice differentiable, so that $t \mapsto
\beta_t([D,T])$ is differentiable, then $[D,T] \in B_{\beta}$, and the main
calculation in the proof of corollary $10.16$ of \cite{GVF} shows that the
derivative of $t \mapsto \beta_t(T)$ at $t=0$ is $[D,T]$.  (We remark 
that it can
be useful, more generally, to use the fact that $\cB(\cH)$ is a von~Neumann
algebra so that the weak operator topology is available.  We can then just ask
that derivatives exist for the weak operator topology.  See in particular
proposition $3.2.53$ of \cite{BrR}.  But we will not pursue this aspect here.)

\begin{lemma}
Suppose that $t \mapsto \beta_t(T)$ is twice differentiable as above.  Then
\[
\|[D,T]\| = \sup\{\|\beta_t(T) - T\|/|t|: t \ne 0\}.
\]
\end{lemma}

\begin{proof}
 From the remarks above it follows that the derivative of $t \mapsto
\beta_t(T)$ is $t \mapsto \beta_t([D,T])$, and that this derivative is
norm-continuous.  Thus
\[
\beta_t(T) - T = \int_0^t \beta_s([D,T])ds,
\]
and so $\|\beta_t(T) - T\| \le |t|\|[D,T]\|$.  This gives inequality in one
direction.  But for any $\epsilon > 0$ we can find $t > 0$ close enough to $0$
that $\|\beta_s([D,T]) - [D,T]\| < \epsilon$ for $0 \le s \le t$. 
Then from the
above integral we obtain
\[
\|\beta_t(T) - T\| \ge t(\|[D,T]\| - \epsilon),
\]
which yields the reverse inequality.
\end{proof}

Again, one can extend the above lemma by using the weak operator topology.  But
for many purposes one wants to deal just with elements which are at least twice
differentiable---see the discussion of {\em regular} tuples in 
section $10.3$ of
\cite{GVF}.

We can now relate our earlier construction of seminorms by means of length
functions to that in terms of ``Dirac'' operators as follows, where the
length function on $\mathbb R$ is the usual one, $\ell(t) = |t|$.

\begin{proposition}
Let $D$ and $B_\beta$ be as above.
Let $A$ be a unital $*$-subalgebra of $B_{\beta}$ such that for any $a \in A$
both $[D,a]$ and $[D,[D,a]]$ are bounded.  Then the seminorms on $A$ defined by
\[
L(a) = \|[D,a]\|
\]
and
\[
L(a) = \sup\{\|\beta_t(a) - a\|/|t|: t \ne 0\}
\]
coincide.
\end{proposition}

In this setting we can thus again ask whether the seminorm $L$ gives $S(A)$ the
weak-$*$ topology.  Usually it will not.  For example, if all of the 
elements of
$A$ commute with $D$ then $L \equiv 0$.  But there are
some quite interesting situations for which it is known
that the answer is affirmative. We give four classes of examples. 

\begin{example}
This first example is the source of the whole
topic, namely the ordinary Dirac operator for a compact spin Riemannian 
manifold $M$,
and Connes' observation that one recovers the ordinary metric on $M$
from the Dirac operator.  More precisely, if $f \in C^{\infty}(M)$ 
and if $f$ is
viewed as an operator by ``pointwise multiplication'' on the spinor bundle with
its Hilbert space structure, then Connes shows that $\|[D,f]\|$ coincides with
the usual Lipschitz seminorm of $f$ from the ordinary metric.  We remark that the
Dirac operator is usually defined for a $\mbox{Spin}^c$ manifold, because these
are the ones which have a spinor bundle.  But if one is only interested in the
smooth structures and ordinary metric of a Riemannian manifold (from which the
Riemannian metric can be recovered), and if the homological 
information which the
Dirac operator contains is not so important, then one can treat any compact
Riemannian manifold.  One simply replaces the spinor bundle (which 
won't exist if
the manifold is not $\mbox{Spin}^c$) with the Clifford algebra bundle itself
equipped with a continuous choice of faithful tracial states so as to give a
Hilbert space structure to the space of cross sections.  Then one uses the
corresponding left regular representation of the Clifford algebra bundle.  The
corresponding Dirac-like operator will again give the smooth structure and the
usual Lipschitz norm for each $f \in C^{\infty}(M)$.
\end{example}

\begin{example}
The type of construction used above
is also used in our second class of examples, which involves the
situation considered right after Theorem 2.2 consisting
of a compact Lie group $G$ acting ergodically on a
unital $C^*$-algebra $A$.  If we put an arbitrary inner-product on the Lie
algebra $\fg$ of $G$, and form the corresponding Clifford algebra,
then the usual construction of a Dirac operator can be
imitated to give a Hilbert space $\cH$, an operator $D$ on it, and a faithful
representation of $A$ on $\cH$.  It is shown in \cite{R4} that the 
corresponding
seminorm on $A$ is a Lip-norm.
\end{example}

\begin{example}
This class of examples concerns the $\theta$-deformed spheres 
and manifolds
of Connes, Landi and Dubois--Violette \cite{CoL} \cite{Cn4} 
\cite{Cn5} \cite{CDV} \cite{Li}. 
These can
be constructed whenever one has an action of the $d$-dimensional 
torus $\bbT^d$,
$d \ge 2$, on any compact manifold $M$.  For any skew-symmetric real $d \x d$
matrix $\theta$ one considers the corresponding non-commutative torus
$A_{\theta}$ with its dual action of $\bbT^d$ which we mentioned earlier.  Then
$\bbT^d$ has a diagonal action on $C(M) \otimes A_{\theta}$, and we let
$M_{\theta}$ be the fixed-point subalgebra for this diagonal action.  This
construction is a reformulation, for the special case of actions of 
$\bbT^d$, of
the general deformation quantization construction for actions of $\bbR^d$
described in \cite{R9}.  This reformulation is also discussed in \cite{Vrl}
\cite{Stz}, and the relation with the quantum groups of \cite{R14} is discussed
in \cite{Vrl}.

Connes and Landi \cite{CoL} \cite{CDV} show that when $M$ is a spin Riemannian
manifold, and when the action $\alpha$ is smooth, leaves the Riemannian metric
invariant, and lifts to the spin bundle, then there is a natural Dirac operator
for $M_{\theta}$.  Hanfeng Li shows in his doctoral thesis \cite{Li} that the
seminorm on $M_{\theta}$ obtained from this Dirac operator is a Lip-norm.
\end{example}

\begin{example}
This class of examples returns to the main example in Connes' 
first paper on
this subject \cite{C1}.  We now let $G$ be a discrete group, and consider its
reduced $C^*$-algebra $C_r^*(G)$ acting on $\boldsymbol{\ell}^2(G)$.  
Let $\ell$ be a length
function on $G$.  We take as our Dirac operator the operator $D = M_{\ell}$ of
pointwise multiplication by $\ell$ on $\boldsymbol{\ell}^2(G)$.  
It is easily seen \cite{C1} that for
any $f \in C_c(G)$, viewed as an operator in $C_r^*(G)$, the operator 
$[D,f]$ is
bounded.  Thus again we can ask if the seminorm $L(f) = \|[D,f]\|$ is a
Lip-norm.  An affirmative answer is now known for two classes of groups.  In 
\cite{R18} it
is shown that if $G = \bbZ^d$ and if $\ell$ is either a word-length function or
the restriction to $G$ of a norm on $\bbR^d$, then the answer is affirmative.
(The proof of this requires a substantial and interesting 
development, involving,
among things, boundaries for non-compact metric spaces.)  In 
\cite{OzR} (and see
also \cite{AnC}) it is shown, by means of techniques which are 
entirely different
from those used in \cite{R18}, that the answer is affirmative if $G$ is a
hyperbolic group and $\ell$ is a word-length function, as well as for certain
reduced free-product $C^*$-algebras.  The techniques involve filtered
$C^*$-algebras and Dirac operators determined by filtrations.  They 
do not apply
to the case of $G = \bbZ^d$.  It would be interesting to find a 
unified proof for
the two cases.
\end{example}

There is ample opportunity to discover additional natural examples of compact
quantum metric spaces.  And there are many further aspects to explore.  For
example, what are the isometry groups \cite{R6} of the above examples?  What 
about quantum
isometry groups \cite{Wng} \cite{Bnc}?  What is the analog of a 
continuous length
function for a compact quantum group, such that it defines Lip-norms on quantum
spaces on which the quantum group acts?

\section{Universality of the Dirac approach}

Although the Dirac operator construction appears fairly special, we 
now show that
in fact every compact quantum metric space can be obtained from the Dirac
operator approach (though this may not be the most useful presentation).
We should make clear now that by this point the term ``Dirac 
operator'' does not
refer to any special kind of self-adjoint operator, but rather to how a
self-adjoint operator is being used, namely to provide metric data by means of
its commutant with the elements of the algebra (or order-unit space) which
specifies the ``space'' being metrized.

We begin by considering ordinary compact metric spaces $(X,\rho)$.  
Let $Z = \{(x,y) \in X \x X: x
\ne y\}$.  Choose any positive measure on $X$ whose support is all of $X$, and
let $\omega$ be the restriction to $Z$ of the square of this measure on $X \x
X$.  Let $A = C(X)$, acting on $L^2(Z,\omega)$ by $(f\xi)(x,y) = 
f(x)\xi(x,y)$.  This is
a faithful representation.  Let $D$ be the operator on 
$L^2(Z,\omega)$ defined by
\[
(D\xi)(x,y) = \xi(y,x)/\rho(y,x),
\]
for those $\xi$ for which $D\xi \in L^2(Z,\omega)$.  Some simple calculations
\cite{R6} show that $[D,f]$ is a bounded operator exactly if $f$ is a Lipschitz
function for $\rho$, and that in this case $\|[D,f]\|$ is exactly the Lipschitz
constant of $f$.  Thus we can recover $\rho$ from $D$.

For the ordinary Dirac operator of a compact Riemannian manifold it is an 
important fact
that it has compact resolvant, that is, the Hilbert space has a basis 
consisting
of eigenvectors, the eigensubspaces have finite dimensions, and any bounded
interval of $\bbR$ contains only a finite number of eigenvalues. 
This led Connes
to require this property of a Dirac operator in his definition of a spectral
triple.  Connes shows that many wonderful properties come from the 
hypothesis of
compact resolvant (such as the existence of a ``volume'' if the spectral triple
is ``$p$-summable'').  In other 
examples the compact resolvant property seems tied to
having some kind of differential structure on $X$.  The operator $D$ 
constructed
earlier for a general compact metric space will usually not have compact
resolvant.  Thus we are led to:

\begin{problem}
Characterize those compact metric spaces for which there is a spectral triple
(for which $D$ has compact resolvant) which gives their metric.  When can the
spectral triple be chosen to be $p$-summable for some $p$?  How does one
characterize the minimal $p$?
\end{problem}

We now turn to the non-commutative case, or more generally the order-unit case.
Let $L$ be a Lip-norm on an order-unit space $A$.  Then $(S(A),\rho_L)$ is an
ordinary compact metric space, and so, as above, there is a representation of
$C(S(A))$ on a Hilbert space $\cH$ and a self-adjoint operator $D$ on 
$\cH$ which
gives $\rho_L$.  But there is a canonical inclusion of $A$ into $C(S(A))$ given
by ${\hat a}(\mu) = \mu(a)$ for $a \in A$ and $\mu \in S(A)$.  (When $A$ is a
$C^*$-algebra, this inclusion is not an algebra homomorphism.)  Under this
inclusion elements of $a$ are carried to functions which are Lipschitz for
$\rho_L$ with Lipschitz constant $L(a)$.  On composing this inclusion with the
representation of $C(S(A))$ on $\cH$ we obtain a faithful representation of $A$
preserving the order unit structure.  Furthermore, $L(a) = \|[D,a]\|$ for every
$a \in A$.

\begin{problem}
For the case of $C^*$-algebras $A$ characterize those Lip-norms which come from
triples $(A,\cH,D)$, where we do not require that $D$ have compact 
resolvant, but
we do require that the representation of $A$ on $\cH$ is a unital $*$-algebra
homomorphism.
\end{problem}

\section{Gromov--Hausdorff distance}

Let $(Z,\rho)$ be an ordinary compact metric space.  For a subset $Y$ 
of $Z$ and
a positive real number $r$ define the open $r$-neighborhood, $\cN_r^{\rho}(Y)$,
of $Y$ by
\[
\cN_r^{\rho}(Y) = \{z \in Z: \mbox{there is a $y \in Y$ with 
$\rho(z,y) < r$}\}.
\]
Hausdorff defined the distance, $\mbox{dist}_H^{\rho}(X,Y)$, between two closed
subsets $X$ and $Y$ of $Z$ by
\[
\mbox{dist}_H^{\rho}(X,Y) = \inf\{r: Y \subseteq \cN_r^{\rho}(X) \mbox{ and } X
\subseteq \cN_r^{\rho}(Y)\}.
\]
This defines an ordinary metric on the set of closed subsets of $Z$, for which
this set is compact.

Gromov generalized this notion of distance to one between any two 
compact metric
spaces. (See \cite{G2}.)  His notion is now called Gromov--Hausdorff distance.
Let $(X,\rho_X)$ and $(Y,\rho_Y)$ be two compact metric spaces.  The
Gromov--Hausdorff distance between them, denoted by $\mbox{dist}_{GH}(X,Y)$, is
defined as follows.  Let $X {\dot \cup} Y$ denote the disjoint union of $X$ and
$Y$.  Let $\cM(\rho_X,\rho_Y)$ be the set of all metrics on the compact set $X
{\dot \cup} Y$ whose restrictions to $X$ and $Y$ are $\rho_X$ and $\rho_Y$
respectively.  Then set
\[
\mbox{dist}_{GH}(X,Y) = \inf\{\mbox{dist}_H^{\rho}(X,Y): \rho \in
\cM(\rho_X,\rho_Y)\}.
\]
For simplicity of notation we are not explicitly indicating the 
metrics $\rho_X$
and $\rho_Y$ on the left-hand side.

There is by now a large and rich literature concerning the
Gromov-Hausdorff limits of compact Riemannian manifolds and related spaces.
As just one example, but one in which the role of Dirac operators is
prominent, and which might suggest phenomena worth investigating in
the non-commutative setting, we mention \cite{Lt}.

We want to define a corresponding notion of Gromov--Hausdorff distance between
compact quantum metric spaces \cite{R6}.  Thus let $(A,L_A)$ and $(B,L_B)$ be
compact quantum metric spaces.  We form the order-unit space $A \oplus B$, and
consider its canonical projections $\pi_A$ and $\pi_B$ onto $A$ and $B$,
respectively.  We let $\cM(L_A,L_B)$ be the set of all Lip-norms $L$ 
on $A \oplus
B$ whose quotient seminorms on $A$ and $B$ are $L_A$ and $L_B$, respectively.
This means that for each $a \in A$ we have
\[
L_A(a) = \inf\{L((a,b)): b \in B\},
\]
and similarly for $L_B$.  There are evident canonical injections of $S(A)$ and
$S(B)$ into $S(A \oplus B)$.  Through these injections we will simply 
view $S(A)$
and $S(B)$ as closed subsets of $S(A \oplus B)$.  From the 
requirement that $L_A$
is the quotient of $L$ it follows (proposition $3.1$ of \cite{R6}) that the
restriction of $\rho_L$ to $S(A)$ coincides with $\rho_{L_A}$, and 
similarly for
$B$.  Thus it is reasonable to define the distance between $(A,L_A)$ and
$(B,L_B)$, which we denote by $\mbox{dist}_q(A,B)$, by
\[
\mbox{dist}_q(A,B) = \inf\{\mbox{dist}_H^{\rho_L}(S(A),S(B)): L \in
\cM(L_A,L_B)\}.
\]
For this to be well-defined and to work smoothly it is 
essential that
the metrics from our Lip-norms give the state spaces the weak-$*$ topology, for
which the state spaces are compact.  This is an important reason for 
our emphasis
in the earlier sections of this paper on proving that the seminorms 
we consider
do give metrics whose topology is the weak-$*$ topology.

An important and non-obvious theorem concerning ordinary Gromov--Hausdorff
distance is that if the Gromov--Hausdorff distance between two compact 
metric spaces is
$0$ then the two metric spaces are isometric \cite{G2}.  Thus Gromov--Hausdorff
distance is actually a metric on the set of isometry-classes of compact metric
spaces.  A similar result holds for our quantum Gromov--Hausdorff distance
(theorem $7.8$ of \cite{R6}), but for this purpose the order-unit spaces of our
compact quantum metric spaces must be suitably completed with respect to their
Lip-norm, so that isometries which should exist will have a place to 
land.  (The
isometries form a compact group in a natural way \cite{R6}.)

However, a defect of the present theory is that two $C^*$-algebras can have
self-adjoint parts which are isomorphic as order-unit spaces, while 
the algebras
themselves are not isomorphic.  This will happen, for example, for those
$C^*$-algebras which are not isomorphic to their opposite algebra.  Thus two
$C^*$-algebras equipped with Lip-norms can be of distance $0$ from each other
without the algebras being isomorphic.

One way of avoiding this defect has been developed by Li \cite{Li}.  He defines
for $C^*$-algebras a variant of our quantum Gromov--Hausdorff distance which
explicitly uses the product in the $C^*$-algebras.  (Thus it has no counterpart
for order-unit spaces.)  Li shows that his variant has many favorable 
properties,
and that it applies to examples of the types which we discuss in the 
next section. Underlying Li's variant is another variant which does apply to
order-unit spaces, and which he calls ``order-unit space quantum Gromov-Hausdorff
distance''. It involves looking at the Gromov-Hausdorff distance between
sets like ${\mathcal L}_1 = \{a \in A: L(a) \leq 1 \ \mathrm{and} \ \|a\|\leq 1\}$
(which are pre-compact, as follows from Theorem 2.1) in the order-unit space
itself, rather than using the state spaces. Li shows that this variant has some 
substantial technical advantages, for example in understanding when a continuous
field of $C^*$-algebras, equipped with Lip-norms in a continuous way, will
be continuous for Gromov-Hausdorff distance.

In section~7 we will discuss another way, due to Kerr \cite{Kr1}, of 
avoiding the
defect.  It involves using operator systems and their matricial norms.

The set of isometry classes of compact quantum metric spaces, equipped with
$\mbox{dist}_q$, is a complete metric space (theorem $12.11$ of \cite{R6}).
There is also a natural analog of a useful theorem of Gromov describing the
totally bounded subsets of this metric space, in terms the number of
$\epsilon$-balls needed to cover $S(A)$ for the various $A$'s in such a subset
(theorem $13.5$ of \cite{R6}).

If $(X,\rho)$ is an ordinary compact metric space, and if $A_X$ is the algebra of
real-valued Lipschitz functions on $X$, then $(A_X,L_{\rho})$ is a 
compact quantum
metric space.  If $(X,\rho_X)$ and $(Y,\rho_Y)$ are compact metric spaces, then
one can show that $\mbox{dist}_q(A_X,A_Y) \le \mbox{dist}_{GH}(X,Y)$.  But
equality may fail.  A simple explicit example of this failure was found by
Hanfeng Li.  (See appendix~1 of \cite{R6}.)  This can be viewed as a defect in
the theory.  But its origin is not difficult to understand.  For ordinary
Gromov--Hausdorff distance one is asking, in effect, that $X$ and $Y$ 
be combined
into a metric space in such a way that the {\em pure} states of $C(X)$ are 
close to the
{\em pure} states of $C(Y)$ (and the other way around); whereas for 
$\mbox{dist}_q$ one
is asking only that each pure state of $X$ be close to some state of $Y$, not
necessarily a pure state.  For non-commutative $C^*$-algebras it is known that
the set of pure states need not be closed, and can be a somewhat bad set, often
fairly inaccessible.  Thus it is not clear how to develop a useful theory in
which one requires pure states of one algebra to be close to pure states of the
other.  On the other hand, for any faithful representation of a 
$C^*$-algebra the
finite convex combinations of vector states will be dense in the 
state space, and
so provide an accessible set of states which can serve well for the theory that
we have developed.  (Note also that for a simple unital infinite dimensional
$C^*$-algebra $A$ the closure of the set of pure states is all of $S(A)$.)

\section{Examples of quantum Gromov-Hausdorff convergence}

For two given ordinary compact metric spaces it is seldom possible to calculate
the Gromov--Hausdorff distance between them precisely.  What has been found
useful in a number of situations is to calculate upper bounds for the distance.
(Lower bounds seem much harder to obtain.)  Then one can try to verify that
certain sequences of compact metric spaces converge to a given space (or form
Cauchy sequences, which must have limit spaces).

The above comments apply equally well for quantum Gromov--Hausdorff distance.
There are now three interesting classes of examples where convergence has been
established.

\begin{example}
For a fixed integer $d$ consider the non-commutative tori 
$A_{\theta}$ discussed
in Section~2, where $\theta$ ranges over real $d \x d$ skew-symmetric matrices.
We have the dual action of $\bbT^d$ on each $A_{\theta}$ as discussed earlier.
Fix a continuous length function on $\bbT^d$, and for each $\theta$ let
$L_{\theta}$ be the corresponding Lip-norm on $A_{\theta}$ according to Theorem
$2.2$.  Then in \cite{R6} it is shown that if a sequence $\{\theta_n\}$ of
matrices converges to a matrix $\theta$, then the sequence
$(A_{\theta_n},L_{\theta_n})$ converges to $(A_{\theta},L_{\theta})$ 
for quantum
Gromov--Hausdorff distance.  In other words, the mapping from matrices $\theta$
to quantum metric spaces $(A_{\theta},L_{\theta})$ is 
continuous for the usual topology on the set of
matrices and for quantum Gromov--Hausdorff distance.
\end{example}

\begin{example}
This example works for any compact semisimple Lie group when suitably rephrased
\cite{R7}, but for simplicity of exposition we describe it only for 
$G = SU(2)$.
Fix a continuous length function on $SU(2)$.  Let $A = C(S^2)$, for $S^2$ the
two-sphere, and consider the usual action of $SU(2)$ on $A$ through the
homomorphism from $SU(2)$ to $SO(3)$ and the action of $SO(3)$ on $S^2$.  This
action on $A$ is ergodic.  Equip $A$ with the corresponding Lip-norm 
as in Theorem $2.2$.  For each $n$ let $(U_n,\cH_n)$ be the irreducible
representation of $SU(2)$ of dimension $n$, and let $A_n = \cB(\cH_n)$, equipped
with the ergodic action consisting of conjugating by $U_n$, as described in
Section~2. (Thus $A_n$ is a full matrix algebra.)  
Then equip each $A_n$ with the corresponding Lip-norm, $L_n$.  In
\cite{R7} it is shown that the sequence $\{(A_n,L_n)\}$ converges to 
$\{(A,L)\}$
for quantum Gromov--Hausdorff distance.  The proof involves Berezin symbols,
which are closely related to coherent states.  This example gives a possible
precise meaning to statements occurring in places in the theoretical physics
literature to the effect that a sequence of matrix algebras converges to the
sphere (or some related space).  See \cite{R7} for references.
\end{example}

\begin{example}
Let the $M_{\theta}$'s be as in Section~4, for a fixed spin Riemannian manifold
$M$ with action of $\bbT^d$, and varying $\theta$'s.  For each $\theta$ let
$L_{\theta}$ be the Lip-norm on $M_{\theta}$ coming from the Dirac operator.
Hanfeng Li shows in \cite{Li} that, much as in Example~$6.1$ above, if
$\{\theta_n\}$ is a sequence of matrices which converges to a matrix $\theta$,
then the sequence $\{(M_{\theta_n},L_{\theta_n})\}$ converges to
$(M_{\theta},L_{\theta})$ for quantum Gromov--Hausdorff distance.
\end{example}

Clearly we are still near the beginning in producing examples of quantum
Gromov--Hausdorff convergence.  But the literature of high-energy
physics and string theory suggests a variety of possible examples.  See the
introduction of \cite{R6} for a number of references.  Furthermore, 
the consequences of
convergence are essentially unexplored in the non-commutative case, but the
discussion of ``degeneration of Riemannian manifolds'' in the literature on
classical Riemannian manifolds (see references in \cite{R6},
and, in particular, \cite{Lt}) suggests many interesting
questions concerning the non-commutative case.

\section{Matricial quantum Gromov--Hausdorff distance}

A matricial version of quantum Gromov--Hausdorff distance has been developed by
David Kerr \cite{Kr1}.  For a variety of reasons it is natural to seek such a
version, but one benefit of it is that it provides one way to repair the defect
of quantum Gromov--Hausdorff distance mentioned in Section~5, namely that two
unital $C^*$-algebras equipped with Lip-norms can have distance $0$ yet not be
isomorphic.

The setting is that of operator systems \cite{ER}.  An operator system is a
self-adjoint subspace of bounded operators on a Hilbert space which 
contains the
identity operator.  (Equivalently, it can be a corresponding subspace of a unital
$C^*$-algebra.)  The self-adjoint part of an operator system is an order-unit
space.  But the crucial difference is that now one is dealing with an 
order-unit
space which has a {\em specific} choice of representation as operators on a
Hilbert space (or in a unital $C^*$-algebra).  The essential feature of this
situation is that if $A$ is an operator system, and if $M_m(A)$ denotes the
linear space of $m \x m$ matrices with entries in $A$ for each $m$, then there is 
a canonical
operator norm and partial order on $M_m(A)$ coming from viewing the matrices as
operators on the sum of $m$ copies of the Hilbert space.  Different (isometric)
representations of an order-unit space can give quite different operator norms
and partial orders on $M_m(A)$ for $m \ge 2$, and thus give distinct operator
systems. 

If $A$ and $B$ are operator systems and if $\varphi$ is a linear map 
from $A$ to
$B$, then we define the corresponding map $M_m(\varphi)$ from $M_m(A)$ to
$M_m(B)$ by $(M_m(\varphi))((a_{ij})) = (\varphi(a_{ij}))$.  We say that
$\varphi$ is {\em $m$-positive} if $M_m(\varphi)$ is positive, and that $\varphi$ is
{\em completely positive} if $M_m(\varphi)$ is positive for all $m$.

Now let $M_n = M_n(\bbC)$ have its canonical matricial structure from its
representation on $\bbC^n$.  We let $UCP_n(A)$ denote the set of unital
completely positive maps from $A$ into $M_n$.  This is the {\em 
$n$-th matricial
state space} of $A$.  It has the evident point-norm topology, for which it is
compact. Notice that $S(A) = UCP_1(A)$.

Now let $L$ be a Lip-norm on (the self-adjoint part of) $A$.  Then for each $n$
we can define a metric, $\rho_{L,n}$, on $UCP_n(A)$ by
\[
\rho_{L,n}(\varphi,\psi) = \sup\{\|\varphi(a) - \psi(a)\|: L(a) \le 1\}.
\]
Kerr shows that because $L$ is a Lip-norm the topology from $\rho_{L,n}$ on
$UCP_n(A)$ coincides with the point-norm topology.  If $\Phi: A 
\rightarrow B$ is
a unital completely positive map from $A$ onto an operator system $B$, then
$UCP_n(B)$ embeds into $UCP_n(A)$ for each $n$ by composition with 
$\Phi$.  Kerr
shows that each of these embeddings is isometric when $UCP_n(B)$ is 
equipped with
the metric from the quotient on $B$ of $L$.

The pieces are then in place to imitate the definition of quantum
Gromov--Hausdorff distance.

\begin{definition}
(definition $3.2$ of \cite{Kr1}.)  Let $(A,L_A)$ and $(B,L_B)$ be Lip-normed
operator systems.  For each $n$ we define the {\em $n$-distance},
$\mbox{dist}_s^n(A,B)$, between $(A,L_A)$ and $(B,L_B)$ by
\[
\mbox{dist}_s^n(A,B) = 
\inf\{\mbox{dist}_H^{\rho_{L,n}}(UCP_n(A),UCP_n(B)): L \in
\cM(L_A,L_B)\}.
\]
If $m > n$ then $\mbox{dist}_s^m(A,B) \ge \mbox{dist}_s^n(A,B)$ (by 
lemma $4.10$
of \cite{Kr1}), and we define the {\em complete distance}, 
$\mbox{dist}_s(A,B)$ by
\[
\mbox{dist}_s(A,B) = \sup_n \{\mbox{dist}_s^n(A,B)\}.
\]
\end{definition}

Kerr shows that these distances satisfy the triangle inequality. 
Furthermore,
he shows that if $\mbox{dist}_s^n(A,B) = 0$ for some $n$, then there is
an isometric $n$-order isomorphism between $A$ and $B$.  For $C^*$-algebras $A$
and $B$ a $2$-order isomorphism will be a $*$-algebra isomorphism (because
by corollary 5 of \cite{Kd2} a unital order isomorphism $\Psi$ will satisfy
$\Psi(a^2) = (\Psi(a))^2$ for self-adjoint $a \in A$, and Choi shows 
(corollary 3.2 of \cite{Cho}) that if $\Psi$ has this latter property and
is also $2$-positive, then $\Psi$ is a *-algebra homomorphism).
Thus:

\begin{theorem}
Let $A$ and $B$ be unital $C^*$-algebras equipped with Lip-norms $L_A$ and
$L_B$.  If $\mbox{\em dist}_s^2(A,B) = 0$ then there is an isometric 
$*$-algebra
isomorphism of $A$ upon $B$ (which carries $L_A$ to $L_B$).
\end{theorem}

Kerr shows that the continuity of quantum Gromov--Hausdorff distance for
non-commutative tori as described in Example $6.1$, and for matrix algebras
converging to the sphere (and other coadjoint orbits) as described in Example
$6.2$, carries over to complete quantum Gromov--Hausdorff distance.

Kerr also provides interesting matricial versions of the theorem on the
completeness of the space of isometry classes of compact quantum metric spaces,
and of the characterization of totally bounded subsets of that space.

\providecommand{\bysame}{\leavevmode\hbox to3em{\hrulefill}\thinspace}
\providecommand{\MR}{\relax\ifhmode\unskip\space\fi MR }
\providecommand{\MRhref}[2]{%
  \href{http://www.ams.org/mathscinet-getitem?mr=#1}{#2}
}
\providecommand{\href}[2]{#2}

\end{document}